\documentclass[10pt,twocolumn,twoside]{IEEEtran}

\usepackage{url}
\usepackage{amsmath}
\usepackage{amsfonts}%
\usepackage{amssymb}%
\usepackage{bm}%
\usepackage{graphicx}
\usepackage{subfigure}
\usepackage{cite}
\usepackage{psfrag}
\usepackage{amsthm}%
\usepackage{mathrsfs}
\usepackage{graphicx}
\usepackage{caption}
\usepackage{mathtools}
\usepackage{color}
\usepackage{algorithm}
\usepackage{algpseudocode}

\theoremstyle{definition}
\newtheorem{theorem}{Theorem}

\newcounter{assumption_}
\newtheorem{assumption}[assumption_]{Assumption}
\newcounter{lemma_}
\newtheorem{lemma}[lemma_]{Lemma}

\newcounter{remark_}
\theoremstyle{remark}
\newtheorem{remark}[remark_]{Remark}

\newenvironment{proof2}[1][Proof]{\textbf{#1.} }{\ \rule{0.5em}{0.5em}}

\title{\Large \bf Event-triggered control cannot improve the $\ell_2$ gain of  $h_\infty$ optimal periodic control and transmit at a smaller average rate}

\author{Duarte J. Antunes and Jo\~ao P. Hespanha
\thanks{Duarte J. Antunes is with the Control Systems Technology Group, Dep. of Mechanical Eng., Eindhoven University of Technology, the Netherlands. {\tt\small d. antunes@tue.nl}. Jo\~ao P. Hespanha is with the Dep. of Electrical and Computer Eng., University of California, Santa Barbara, CA 93106-9560, USA. {\tt\small hespanha@ece.ucsb.edu}.  }} 

\begin{document}

\maketitle
\thispagestyle{empty}
\pagestyle{empty}

\begin{abstract}
We consider a standard discrete-time event-triggered control setting by which a scheduler collocated with the plant's sensors decides when to transmit sensor data to a remote controller collocated with the plant's actuators.  When the scheduler transmits periodically with period larger than or equal to one, the $h_\infty$ optimal controller guarantees an optimal attenuation bound ($\ell_2$ gain) from any square-summable disturbance input  to a plant's output. We show that, under mild assumptions, there does not exist a controller and scheduler pair  that strictly improves the optimal attenuation bound of periodic control with a smaller average transmission rate. Equivalently, given any controller and scheduler pair, there exists a square-summable disturbance such that either the attenuation bound or the average transmission rate are larger than or equal to those of optimal periodic control. 
\end{abstract}
\section{Introduction}
\par Event-triggered control (ETC) has been the subject of extensive research over the last two decades~\cite{astrom:02,xu:04,Tabuada2007ETC,molin:10,lunze:11,ramesh:11,heemels:12_2,postoyan:15,chen:12,antunes2014,araujo2014cdc,antunes:16_1,behnam:18,antunes:20,hadi:18,gatsis:2016,mirkin:17,victor:16,hadi:20,mirkin:19}. It provides an alternative to traditional periodic control that can potentially reduce the communication burden for the same closed-loop performance or, equivalently, improve the closed-loop performance using the same communication resources. This desired feature has been addressed in the literature considering two measures of performance inherited from the standard $h_\infty$ and $h_2$ frameworks for periodic linear systems~\cite{chen:95}.
\par In the $h_2$ control setting, closely related to Linear Quadratic Gaussian Control (LQG), performance is measured by an average quadratic cost. Considered in this setting, ETC can (strictly) improve the  average quadratic cost of optimal periodic control for the same, and even smaller, transmission rate;  see~\cite{astrom:02} for scalar  systems and ~\cite{chen:12,antunes2014,antunes:16_1,araujo2014cdc,gatsis:2016,hadi:18,behnam:18,antunes:20} for general systems. This property is often referred to as consistency~\cite{antunes:16_1} or strict consistency if the performance can be~\textit{strictly} improved~\cite{antunes:20}.
\par In the $h_\infty$ setting, disturbances are assumed to be deterministic and performance is defined through the $\ell_2$ gain. This gain captures the smallest attenuation bound from the disturbance input to an output of interest of the plant and coincides with the $h_\infty$ norm when the system is linear. While much work on ETC has considered the $\ell_2$ gain~\cite{heemels:12_2,postoyan:15}, obtaining consistency properties in this setting similar to the ones obtained in the $h_2$ setting has received little attention. Two exceptions are~\cite{hadi:20} and~\cite{mirkin:19};~\cite{mirkin:19} considers continuous-time linear systems and an approach based on Youla parametrizations, whereas~\cite{hadi:20}  considers discrete-time systems and tackles this problem with a game theoretical approach. Given a desired attenuation bound that is guaranteed by periodic control with a given transmission rate,  they provide an ETC policy that can guarantee the same attenuation bound with a smaller or equal average transmission rate. Still, one can select the disturbance input so that transmissions are triggered periodically, resulting in the same $\ell_2$ gain. In this sense, the event-triggered controllers in~\cite{hadi:20},~\cite{mirkin:19} are not strictly consistent.
\par Motivated by this recent research, in this paper, we pose the following question: given a discrete-time linear system and optimal periodic controller with transmission rate between sensors and actuators $\frac{1}{h}$, for some period $h \in \mathbb{N},$ and $h_\infty$ norm $\gamma_h$ (smallest attenuation bound), can one find  scheduling and controller policies that lead to a strictly smaller attenuation bound than $\gamma_h$ while transmitting at a smaller average rate $\frac{1}{h}$. The answer is no, in the sense that, given any controller and scheduling policies and period $h \in \mathbb{N},$, there exists a disturbance such that either the attenuation bound is larger than or equal to $\gamma_h$ or the transmission rate is smaller than or equal to $\frac{1}{h}$. This result is illustrated in Figure~\ref{fig:2} below. The proof of this result is constructive in the sense that we provide a disturbance  policy that generates such a disturbance. This is a fundamental limitation of ETC in the context of $h_\infty$ control, which is not present in $h_2$/LQG control. 
\par Note that for a \textit{given} disturbance found in a practical setting, a given controller-scheduler pair (such as the one proposed in~\cite{hadi:20}) can result in both smaller attenuation from \textit{that disturbance} to a plant's output and smaller average transmission rate when compared to  the optimal periodic $h_\infty$ controller with rate $\frac{1}{h}$, as desired. From the perspective of our result, this desired property cannot hold for \textit{every} disturbance.
\par The paper is organized as follows. Section~\ref{sec:2} formulates the problem and Section~\ref{sec:4} provides the main results. Section~\ref{sec:5} provides numerical examples and Section~\ref{sec:6} gives concluding remarks. The proofs of the results are given in Section~\ref{sec:7}.
\section{Problem formulation}\label{sec:2}
\par Consider a linear system 
\begin{equation}\label{eq:sys}
	x_{t+1} = Ax_t+B_2u_t+B_1w_t, \ \ t\in \mathbb{N}_0:=\mathbb{N} \cup \{0\},
\end{equation}
with the following output of interest
\begin{equation}\label{eq:output}
z_{t} = C_2 x_t+D_{21}u_t,
\end{equation}
where $x_t \in \mathbb{R}^n$, $u_t \in \mathbb{R}^m$, $z_t \in \mathbb{R}^p$, $w_t \in \mathbb{R}^{n_w}$ for $t\in \mathbb{N}_0$. Without loss of generality, we assume that  $C_2^\top D_{21}=0$, and define $Q:=C_2^\top C_2$ and $R = D_{21}^\top D_{21}$. Futhermore, we consider $B_1=I$ and use the notation $B=B_2$. The following standard assumption is stated for future reference.
\begin{assumption}\label{as:1}
	 $(A,B)$ is controllable, $(A,C_2)$ is observable and $R>0$. \hfill $\square$ 
\end{assumption} Let  $w=(w_0,w_1,w_2,\dots)$, $z=(z_0,z_1,z_2,\dots)$,  define the inner product $\langle w,z\rangle=\sum_{t=0}^\infty w_t^\top  z_t$ and norm
$\|w\|^2:=\sqrt{\langle w,w\rangle}$, and let $\ell_2$ be the Hilbert space of sequences with bounded norm. The system provides an attenuation bound $\gamma$ from the input disturbances to the output of interest if 
\begin{equation}\label{eq:inequality}
	\|z\|^2\leq \gamma^2\|w\|^2, \ \ \forall w \in \ell_2, \text{ for } x_0=0.\ \ \  
\end{equation}
We are interested in ensuring that~\eqref{eq:inequality} holds for the smallest possible $\gamma$.  We assume $w_0 \neq 0$ so that $x_1=w_0$; there is no loss of generality when the control policy sets $u_t=0$ when $x_t=0$, which is typically the case. Then, $w_t=0$ for all $t$ leads to $z_t=0$ for all $t$ and~\eqref{eq:inequality} is trivially met; otherwise, the first non-zero component of $w$ can be assumed to be $w_0$, due to time invariance of~\eqref{eq:sys}. The initial condition $x_0$ may be non-zero provided that we redefine~\eqref{eq:inequality} as, e.g., in~\cite[Sec.~3.5.2]{basar:91}.

\par
\begin{figure}
	\centering
	\includegraphics[width=6cm]{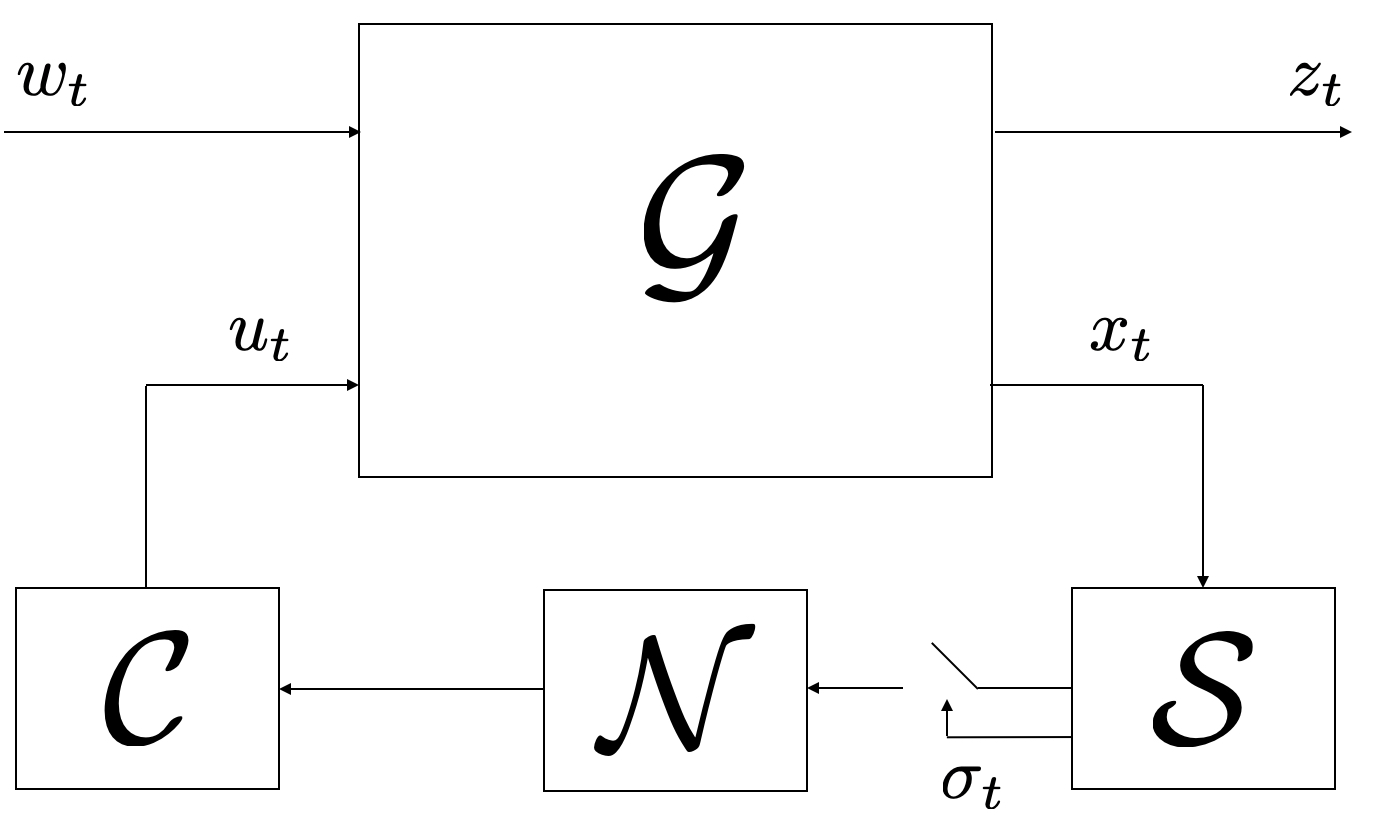}
	\caption{An event-triggered controller consists of a controller $\mathcal{C}$ and a scheduler $\mathcal{S}$, which sends measurement/state data to the controller; $\mathcal{N}$ represents the network and $\mathcal{G}$  the plant.}\label{fig:1}
\end{figure}

 We consider the networked control system depicted in Figure~\ref{fig:1}. The control input $u_t$ is determined by two agents:  the scheduler which controls the information seen by the controller since it decides when to transmit sensor/state data to the controller over a network, and the controller which provides $u_t$ at every $t$. Let $\sigma_t = 1$ if the scheduler transmits the state to the controller at time $t$ and $\sigma_t =0$ otherwise; at $t=1$ it is assumed that the scheduler sets $\sigma_1=1$. Moreover, let  $s_\ell$ be the transmission times defined as $s_{\ell+1}=s_\ell+\tau_\ell$, $s_0=1$ with $\tau_\ell = \min\{j \in \{1,\dots,\bar{h}\}|\sigma_{s_\ell+j}=1\}$  where $\bar{h}$ is a given integer implicitly imposing the following assumption.
\begin{assumption}\label{as:2}
	The scheduler is such that $\tau_\ell \leq \bar{h} $ for a given $\bar{h} \in \mathbb{N}$. \hfill $\square$
\end{assumption}
 This is a mild assumption since $\bar{h}$ can be arbitrarily large, and any scheduler can be modified to trigger when $r_t = \bar{h}$, where $r_t = t-s_{\bar{\ell}(t)}$ is the elapsed time since the last transmission, $s_{\bar{\ell}(t)}$ is the time of the most recent transmission prior to time $t$ and ${\bar{\ell}}(t)=\max\{\ell \in \mathbb{N}_0|s_\ell \leq t|\}$. 
 Let $ \mathcal{I}_t:= \{x_k| k \in \{s_{\bar{\ell}(t)},s_{\bar{\ell}(t)+1}\dots,t\} \}$ be information sets containing the states from the previous transmission time up to the current time $t$  and  $ \mathcal{J}_t:= \{x_{s_{\bar{\ell}(t)}}\}\cup\{r_t\}$ be information sets containing only the last transmitted state and $r_t$.  The information available to the scheduler, controller and disturbance policies is summarized in the next assumption.
  \begin{assumption}\label{as:3}
  	The scheduler, controller and disturbance policies are assumed to depend on the following information sets for every $t \in \mathbb{N}_0$ 
  	\begin{itemize}
  	 \item[(i)] $ \sigma_t = \mu_{\sigma,t}(\mathcal{I}_t)$,
  	 \item[(ii)] $ u_t = \mu_{u,t}(\mathcal{J}_t)$,
  	 \item[(iii)] $ w_t = \mu_{w,t}(\mathcal{I}_t).$  \hfill $\square$
  	\end{itemize}
  \end{assumption}
  While it is typical to assume that the information sets include states of the system since the initial time $t=0$ and then show that the obtained policies rely only on recent states (see, e.g.,~\cite{hadi:20},~\cite[Ch.~3]{basar:91}), here we immediately assume the  policies can only depend on most recent relevant states. Thus, this is also a mild assumption. Its main purpose is to formulate Assumption~\ref{as:5} below. Two general examples of controller and scheduler pairs are given in Section~\ref{sec:5} and Algorithm~\ref{alg:1} is an example of a disturbance generator policy. We assume  the controller policy is known to the scheduler and to the disturbance policy, and thus it would be redundant to add $u_k$, $k \in \{s_{\bar{\ell}(t)},\dots,t\}$, to $\mathcal{I}_t$. Adding $r_t$ to $\mathcal{I}_t$ would also be redundant as it is inferred by the number of states in $\mathcal{I}_t$.
  
  \par We define $\pi_a = (\mu_{a,0}, \mu_{a,1},\dots)$ as the policy of the event-triggered controller when $a=\sigma$, the policy of the controller when $a=u$ and the policy of the disturbance generator when $a=w$. Given  $\pi := (\pi_u,\pi_\sigma)$, we define the average transmission rate 
$ r_{\pi}(w):=\limsup_{s\rightarrow \infty} \frac{1}{s}\sum_{t=0}^{s-1}\sigma_t$
and the average inter-transmission interval $ \bar{h}_{\pi}(w) = \frac{1}{ r_{\pi}(w)}$.

\subsection{Periodic control}\label{subsec:B}
\noindent Periodic schedulers with evenly spaced samples are defined as
\begin{equation}\label{eq:sch}
\sigma_t=\left\{\begin{aligned}
&1 \text{ if }t \text{ is zero or an integer multiple of }h,\\
&0 \text{ otherwise. }
\end{aligned}\right.
\end{equation}
As explained in Remark~\ref{rem:2} below,~\eqref{eq:sch} are superior in a given sense to other periodic schedulers and we can restrict our attention to~\eqref{eq:sch} for comparison with event-triggered schedulers. For brevity, we denote~\eqref{eq:sch} by periodic schedulers.
\par Note that for this scheduler $r_{\pi}(w)=\frac{1}{h}$ for every $w$. Let 
\begin{equation}\label{eq:gammah}	\begin{aligned}
		\gamma_h:= \inf\{\gamma| \exists \pi_u & \text{ such that }\eqref{eq:inequality}\text{ holds} \\ & \text{ when the scheduler is given by~\eqref{eq:sch}}\}.
	\end{aligned} 
\end{equation}
  To compute $	\gamma_h$ let us first define three matrix transformations:
$$ \begin{aligned}
F_a(P)&:= P+P(\gamma^2 I-P)^{-1}P \\
F_c(P)&:= A^\top PA+Q-A^\top PB(B^\top PB+R)^{-1}B^\top PA  \\
F_o(P)&:= A^\top PA+Q.
\end{aligned}$$
 for a given $\gamma \in \mathbb{R}_{>0} $.  When $h=1$, and under Assumption~\ref{as:1}, the following iteration 	$P_{t+1}=F_c(F_a(P_t))$ 
with $P_0=0$ is monotone in the sense that $P_{t+1}\geq P_t$ converges if $\gamma^2 I>P_t$ for every $t \in \mathbb{N}_0$ to the unique positive definite solution $\bar{P}_\gamma$ of the algebraic Riccati equation 
\begin{equation}\label{eq:Pgamma} 	\bar{P}_\gamma=F_c(F_a(\bar{P}_\gamma)), \end{equation}
see~\cite[Sec~3.2]{basar:91} (although the expressions in~\cite[Sec~3.2]{basar:91} appear in a different but equivalent form). Due to monotonicity, $\gamma^2 I>\bar{P}_\gamma$ implies that  $\gamma^2 I>P_t$ for every $t \in \mathbb{N}_0$. Provided that this condition holds,~\eqref{eq:inequality} holds for a policy $\pi_x$ specified by $u_t= K_\gamma x_t$, 
where 
\begin{equation}\label{eq:K}
	K_\gamma= -(R+B^\top F_a(\bar{P}_\gamma)B)^{-1}B^\top F_a(\bar{P}_\gamma)A.
\end{equation}
%
If $\gamma$ is such that $\gamma^2 I \geq P_{t}$ does not hold for some $t$, then~\eqref{eq:inequality} does not hold for any $\pi_u$.
\par However, if $h>1$ the conditions on $\gamma$ for the existence of such a control policy are stricter, i.e., $\gamma$ needs to be larger~\cite{hadi:20}. They actually become stricter as $h$ increases leading to non-decreasing sequence of $\gamma_h$. In fact,
 consider the following iteration, with $M_1=\bar{P}_\gamma$,
	\begin{equation}\label{eq:iter}
		M_{k+1} = F_o(F_a(M_k)), \ \ \ k \in \{1,2,\dots,{h}\},
	\end{equation} 
	with $\gamma$ such that $\gamma^2 I\!-\!M_{k}\!>\!0$, for  all $k \!\in\! \{1,\dots,h\}$. Then, as proven in~\cite[Lemma 1]{antunes:2024}, $M_{k+1}\!\geq\! M_k$,  $\inf\{ \gamma|\gamma^2 I-M_h>0\}=\gamma_h$, with $\gamma_h$ given in~\eqref{eq:gammah}, and $\gamma_{h+1}\geq \gamma_h$ for all $h \in \mathbb{N}$.

\subsection{Problem statement}\label{subsec:C}
\begin{figure}[t!]
	\centering
	\includegraphics[width=8.5cm]{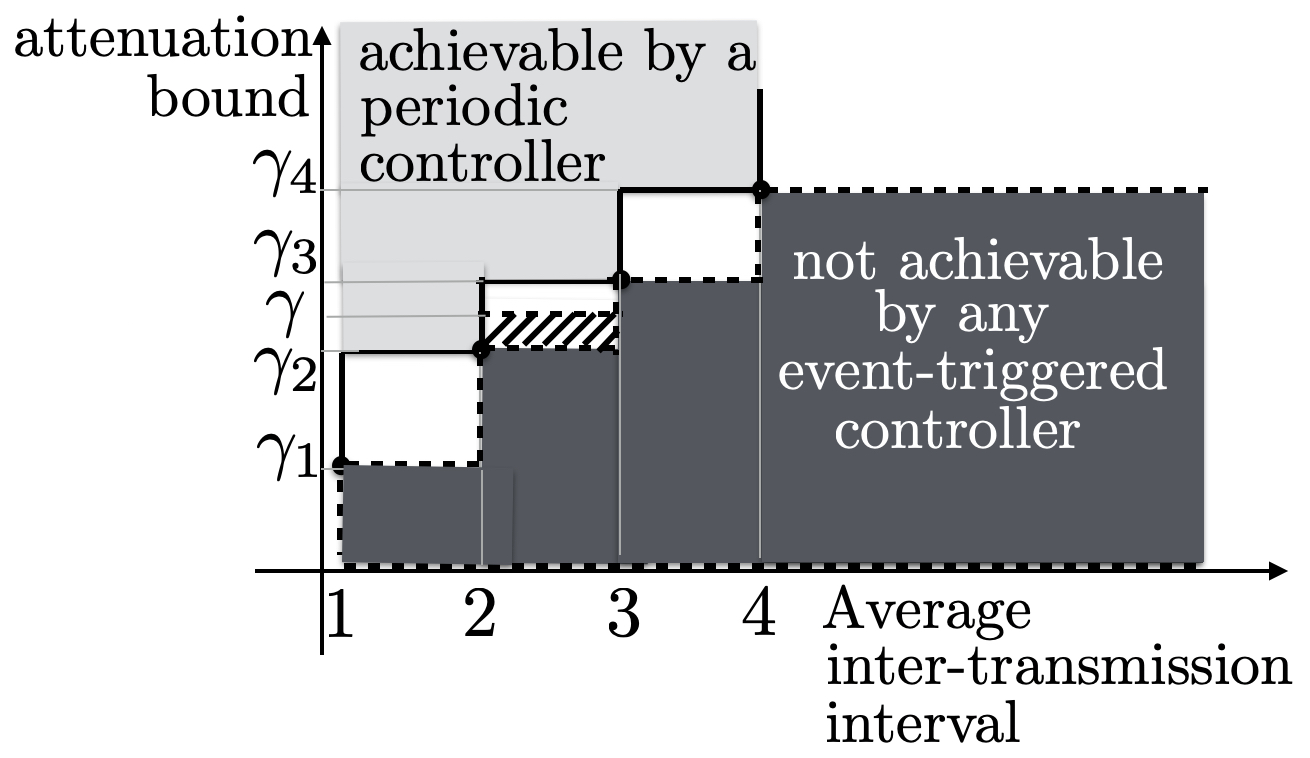}
	\caption{Illustration of our main result Theorem~\ref{th:1}; $\gamma_h$ is the $\ell_2$ gain (smallest attenuation bound which coincides with the $h_\infty$ norm) of periodic control with period $h$. See  Remarks~\ref{rem:1},~\ref{rem:2}.}\label{fig:2}
\end{figure}
\par For a given control and scheduling policy $\pi=(\pi_u,\pi_\sigma)$ let
\begin{equation}\label{eq:gammapibar}
	 \begin{aligned}
	&\!\!\!\! \gamma_p \! :=\!\inf\{\gamma |\text{ such that}~\eqref{eq:inequality} \text{ holds  when } u_t = \mu_{u,t}(\mathcal{J}_t), \\
	&\ \ \ \ \ \ \ \ \ \ \ \ \ \ \ \ \ \ \ \  \hspace{3cm} \sigma_t = \mu_{\sigma,t}(\mathcal{I}_t)\}.
	 \end{aligned}
\end{equation}
In this paper we pose the question of whether, for a given~$h$,  it is possible to find policies $\pi_u$ and $\pi_\sigma$ such that 
 \begin{itemize}
	\item[(i)] $\gamma_\pi<\gamma_h$,  and 
	\item[(ii)]  $r_{\pi}(w) <\frac{1}{h}$, for every $w \in \ell_2$.
\end{itemize}
\par Graphically, the question is whether we can find a controller and scheduler pair such that the pair $(r_\pi,\gamma_\pi)$ is in the dark gray region of Figure~\ref{fig:2}, which is an open set. The answer turns out to be no as stated in Theorem~\ref{th:1} below. 
\begin{remark}\label{rem:1}
			Theorem~\ref{th:1} is actually stronger in the sense that if we can find $\gamma$ such that $\gamma_h<\gamma<\gamma_{h+1}$ and  Assumptions~\ref{as:4}-\ref{as:5} are satisfied, then the pair $(r_\pi,\gamma_\pi)$ cannot additionally lie in the region with diagonal stripes. If Assumptions~\ref{as:4}-\ref{as:5} are satisfied for every $h$ and $\gamma$ such that $\gamma_h<\gamma<\gamma_{h+1}$ the pair $(r_\pi,\gamma_\pi)$ of event-triggered control cannot  lie in the union between the dark gray and white regions, again an open set.
	\end{remark}
\begin{remark}\label{rem:2}
	In~\cite{antunes:2024} it is proven that the $\ell_2$ gain associated with a periodic scheduler with largest interval between samplings/transmissions $\tilde{h}$ is equal to the $\ell_2$ gain associated with~\eqref{eq:sch} with $h=\tilde{h}$.  Thus,  a periodic scheduler with non-evenly spaced sampling is inferior to~\eqref{eq:sch}  as it  results in an average inter-transmission interval  smaller than $h$ and $\ell_2$ gain $\gamma_h$. This latter fact is also implied by Theorem~\ref{th:1} below since we allow for time-varying scheduling policies. Given a desired rational (dense in the reals) average inter-transmission interval smaller than $h$  (e.g., $9/4$ with $h=3$), we can always design a periodic scheduler with that average inter-transmission interval and with $\ell_2$ gain $\gamma_h$ (e.g., repeat schedules 100101010). Thus, we can guarantee any average inter-transmission interval and attenuation bounds in the interior of the light gray region with periodic control. 
\end{remark}


\section{Main result}\label{sec:4}
  

\par The main result of the paper, stated in this section, compares ETC to periodic control with an arbitrary but fixed period $h \in \mathbb{N}$. Besides Assumptions~\ref{as:1}-\ref{as:3}, it needs the following two additional technical assumptions on the given period $h$.
\begin{assumption}\label{as:4}
	$h$ is such that $\gamma_h< \gamma_{h+1}$. \hfill $\square$
\end{assumption} 
\par  Note that Assumption~\ref{as:4} can easily be tested. From ~\cite[Lemma 1]{antunes:2024} we know that $\gamma_h\leq \gamma_{h+1}$, but we require a strict inequality to state Assumption~\ref{as:5} below.  While Assumption~\ref{as:4} typically holds one can find cases where it does not hold. For example if $A=0$ then $F_o(P)=Q$, $M_k=Q$ for every $k$ in~\eqref{eq:iter}, and $\gamma_h=\gamma_{h+1}$ for every $h$. 
\par When Assumption~\ref{as:4} holds, consider $\gamma$ such that $\gamma_h< \gamma<\gamma_{h+1}$. Then $(\gamma^2I-\bar{P}_\gamma)$ is invertible and we can define  \begin{equation}\label{eq:Lgamma}
	L_{\gamma}:=(\gamma^2I-\bar{P}_\gamma)^{-1}\bar{P}_\gamma A. \end{equation}
Moreover,  the smallest eigenvalue of $\gamma^2I-M_{h+1}$ is negative, and we denote by
$v_{h,\gamma}$  (either) one of the two unitary-euclidean norm eigenvectors associated with this smallest eigenvalue. For each transmission time $t=s_\ell$ for some $\ell$, and for a given $\epsilon \geq 0$  consider the following disturbance policy to be used until the next transmission time $s_{\ell+1}$,
\begin{equation}\label{eq:w2}
	w_k \!=\!	\left\{\begin{aligned} &\! L_\gamma (Ax_k\!+\!Bu_k)+\epsilon v_{h,\gamma}, \text{ if }k = s_\ell\\
	&\!	L_\gamma (Ax_k\!+\!Bu_k), \text{if } \tau_\ell \!>\!1 \!\text{ and }\! k \!\in\!  \{s_\ell\!+\!1,\dots,s_{\ell\!+\!1}\!-\!1\}.
	\end{aligned}\right.
\end{equation}
When $\epsilon=0$ this is a well-known optimal disturbance policy for a game between control and disturbances with payoff $ \|z\|^2 -\gamma \|w\|^2$. Let $\mathcal{T}(\xi,t,\epsilon,\gamma)=\tau_\ell$ denote the inter-transmission time when~\eqref{eq:w2} is applied at time $t=s_\ell$ and $x_{s_\ell}=\xi$, which only depends on the state $\xi$ and time $t$ since the control policy is assumed to be known. The dependency on the disturbance parameters $\epsilon$, $\gamma$ is added for convenience. We need the following regularity property on this map.
\begin{assumption}\label{as:5}
	Under Assumption~\ref{as:4}, there exist $\gamma$, with $\gamma_h< \gamma <\gamma_{h+1}$, and $\underline{\epsilon}>0$,  such that for every $\xi \in \mathbb{R}^n$ and $t \in \mathbb{N}$  either 
		\begin{equation}\label{eq:A}
		\mathcal{T}(\xi,t,0,\gamma)\leq h
	\end{equation}
	or
	\begin{flushright}
	$\mathcal{T}(\xi,t,\epsilon,\gamma)> h, \text{ for every }\epsilon \in [-\underline{\epsilon},\underline{\epsilon}]$. \ \ \ \ \ \ \ \ \ \   $\square$
	\end{flushright} 
\end{assumption} 

This assumption states that, for every state $x_t=\xi$  and every transmission time $t=s_\ell$ either the scheduler triggers/transmits before (and including) time $s_{\ell}+h$ when the disturbance policy~\eqref{eq:w2} with $\epsilon=0$ is used or, otherwise, does not trigger before (and including) time $s_{\ell}+h$ when~\eqref{eq:w2} is used for an arbitrarily small $\epsilon$ (besides $\epsilon=0$). However, $|\epsilon|$ must be uniformly lower bounded with respect to $\xi$ and $t$. 
\par The numerical examples illustrate how to test Assumption~\ref{as:5}.  Even for a general scheduler and controller pair, this assumption is mild in the sense that, given any scheduler and controller pair that does not meet Assumption~\ref{as:5}, we can apply a small modification to the scheduler to ensure that Assumption~\ref{as:5} is met. To this effect it suffices for each $\xi$ and $t$ for  which~\eqref{eq:A} does not hold to set the scheduler to $\mu_{\sigma,t}(\mathcal{I}_t)=0$ for the state trajectories in $\mathcal{I}_t$ resulting from~\eqref{eq:w2} with $\epsilon \in [-\underline{\epsilon},\underline{\epsilon}]$; this is a very small set when compared to its complement and commensurable with $\underline{\epsilon}$.  

\par  We are ready to state the main result. 

\begin{theorem}\label{th:1}
	Consider linear system~\eqref{eq:sys} with performance output~\eqref{eq:output}, an arbitrary $h\in \mathbb{N}$ and suppose that Assumptions~\ref{as:1}-\ref{as:5} hold for some $\gamma$ such that $\gamma_h< \gamma < \gamma_{h+1}$. 
  	  Then,  there exists $w\in \ell_2$ such at least one of the following two must hold:
	\begin{itemize}
		\item[(i)]   $\|z\|^2>  \gamma^2\|w\|^2$, which implies that $ \gamma_\pi \geq \gamma> \gamma_h$, 
		\item[(ii)] $r_{\pi}(w) \geq  \frac{1}{h}.$ \hfill $\square$
	\end{itemize}
\end{theorem}

\begin{proof2} Algorithm~\ref{alg:1}, associated with a $\gamma$ for which Assumptions~\ref{as:4},~\ref{as:5} hold and initialized with an arbitrary $w_0$, provides a disturbance policy that when applied to~\eqref{eq:sys} with a given scheduler and controller pair, fixed disturbance values $w_1, w_2, \dots$ are generated such that either $ \gamma_\pi > \gamma> \gamma_h$
 or $r_{\pi}(w) \geq  \frac{1}{h}$. 
Letting $\eta_t\! =\! \sum_{j=0}^{t-1}\!z_j^\top z_j\!-\!\gamma^2w_j^\top w_j \!+\!x_t^\top \bar{P}_{{\gamma}}x_t$, if
\begin{equation}\label{eq:ineq} \sum_{j=0}^{\infty}z_j^\top z_j-\gamma^2w_j^\top w_j=\eta_t-x_t^\top \bar{P}_{{\gamma}}x_t+\sum_{j=t}^{\infty}z_j^\top z_j-\gamma^2w_j^\top w_j>0
\end{equation}
the $\ell_2$ gain cannot be strictly smaller than $\gamma$, i.e., $ \gamma_\pi \geq \gamma> \gamma_h$.   Algorithm~\ref{alg:1}, at each transmission time $t=s_\ell$ for which $\eta_t\leq 0$, checks~\eqref{eq:A} by probing the system with~\eqref{eq:w2} with $\epsilon =0$ (this can be done causally through simulation since the model, the control and scheduling policy are known). This is an optimal disturbance policy in the sense mentioned below~\eqref{eq:w2}. If~\eqref{eq:A} holds, then Algorithm~\ref{alg:1} does apply~\eqref{eq:w2} with $\epsilon =0$, i.e.,~\eqref{eq:wA1}. However, when~\eqref{eq:A} does not hold,~\eqref{eq:wA1} is not necessarily optimal for the disturbance generator. In fact, the disturbance policy can increase $\sum_{j=t}^{t+\tau_\ell}z_j^\top z_j-\gamma^2w_j^\top w_j$ (which contributes to increasing $\sum_{j=t}^{\infty}z_j^\top z_j-\gamma^2w_j^\top w_j$ in~\eqref{eq:ineq}) when compared to the case where~\eqref{eq:wA1} would be used for every $\ell$, by applying~\eqref{eq:w2} with  $\epsilon \neq 0$,\textit{ provided that the control inputs} $u_{s_\ell},\dots,u_{s_\ell+h}$ \textit{ do not change} and still guarantee $\tau_\ell>h$. This follows from  Lemma~\ref{lem:key} below. Algorithm~\ref{alg:1} computes a set of $\epsilon$ for which $\tau_\ell>h$ which implies $u_{s_\ell},\dots,u_{s_\ell+h}$ do not change. Assumption~\ref{as:5} ensures that this set contains at least $\epsilon \in [-\underline{\epsilon},\underline{\epsilon}]$. Thus,  if $\tau_\ell > h$ then~\eqref{eq:w2} with $\epsilon \neq 0$, i.e.,~\eqref{eq:wA2}, is used where $\epsilon$ is computed to ensure that $\sum_{j=t}^{t+\tau_\ell}z_j^\top z_j-\gamma^2w_j^\top w_j$ is increased when compared to the case where~\eqref{eq:w2} with $\epsilon=0$ would be applied. This is the rationale behind the expression~\eqref{eq:eps} for $\epsilon$ which indeed increases this cost taking into account~\eqref{eq:FF1}. This cost increase contributes to increasing $\eta_{s_{\ell+1}}$. Algorithm~\ref{alg:1} performs these steps while monitoring if 
$\eta_t>0$.  When $\eta_t>0$,~\eqref{eq:wopt_} is used. Then, relying on Lemma~\ref{lem:2} below	$x_t^\top \bar{P}_\gamma  x_t -\eta_t/2 \leq  \sum_{j=t}^{\infty}z_j^\top z_j-\gamma^2w_j^\top w_j$. 
Combining this inequality and $\eta_t>0$  we conclude that
\begin{equation}\label{eq:ineqtruncated2}
	\underbrace{\sum_{j=0}^{t-1}z_j^\top z_j-\gamma^2w_j^\top w_j}_{=\eta_t -x_t^\top \bar{P}_\gamma x_t} +	\underbrace{\sum_{j=t}^{\infty}z_j^\top z_j-\gamma^2w_j^\top w_j }_{\geq x_t^\top \bar{P}_\gamma x_t-\eta_t/2}\geq \eta_t/2>0
\end{equation}
\noindent which is~\eqref{eq:ineq}. It is clear that $w\in\ell_2$ since $w_t=0$ after a finite time. If $\eta_t\leq 0$ for every $t=s_\ell$, $\ell \in \mathbb{N}$, Lemma~\ref{lem:final} shows that $r_{\pi}(w) \geq  \frac{1}{h}$ and  $\omega \in \ell_2$, concluding the proof.
\end{proof2}
\begin{algorithm}[t!]
	\caption{Policy of disturbance generator}\label{alg:1}
	\begin{algorithmic}[1]
		\State Choose an arbitrary $w_0 \neq 0$, an arbitrary $\bar{\epsilon}\geq \underline{\epsilon}> 0$, and compute $v_h$ and $\gamma$ such that Assumption~\ref{as:5} holds.
		\For{ $\ell \in \{0,1,2,\dots\}$ (at time $t=s_\ell$, $s_0=1$)}
		\State Compute  $\eta_t = \sum_{j=0}^{t-1}z_j^\top z_j-\gamma^2w_j^\top w_j +x_t^\top \bar{P}_{{\gamma}}x_t$
		\If{$\eta_t>0$} \textbf{break (go to line 11)}
		\Else 	\If{ ~\eqref{eq:A} holds  with $\xi=x_{s_{\ell}}$, $t=s_{\ell}$} 
		set
		\begin{equation}\label{eq:wA1}
			\begin{aligned}
				\ \ \ \ \ \ \ \ \ \ \ w_k=\bar{L}_\gamma(Ax_k+Bu_k),\text{ if }k\in \{t,\dots,s_{\ell+1}\} 
			\end{aligned}
		\end{equation}	
	\ \ \ \ \ \ \ \ \ \ \ \ 	(then $s_{\ell+1}=t+\mathcal{T}(x_{s_{\ell}},s_{\ell},0,\gamma)$)
		\Else   { Compute} $\mathcal{E}\!=\!\{\epsilon \!\in\! \mathbb{R}|\mathcal{T}(\xi,t,\epsilon,\gamma)\!>\! h\}\!\cap\! [-\bar{\epsilon},\bar{\epsilon}]$,
		\begin{equation}\label{eq:eps}
			\epsilon =\left\{\begin{aligned}
				& \sup \mathcal{E} \text{ if }b(x_{s_\ell},U_{s_\ell})\geq 0\\
				& \inf \mathcal{E} \text{ if }b(x_{s_\ell},U_{s_\ell})< 0
			\end{aligned}\right.
		\end{equation}
	\ \ \ \ \ \ \ \ \ \ \ \ \ \ 	where $b(x_t,U_t)$ is given by~\eqref{eq:FFFF2} below and set
		\begin{equation}\label{eq:wA2}
		\ \ \ \ \ \ \ \	w_k=\left\{\begin{aligned}
				&	\bar{L}_\gamma(Ax_k+Bu_k)+\epsilon v_{h,\gamma},\text{ if }k =t, \\
				&	\bar{L}_\gamma(Ax_k+Bu_k),\text{ if }k\in \{t+1,\dots,s_{\ell+1}\} 
			\end{aligned}\right.
		\end{equation}	
		\EndIf
		\EndIf
		\EndFor
		\State Compute $q=\zeta(x_k,\eta_k/2)$ with the method in Lemma~\ref{lem:2} and set 
		\begin{equation}\label{eq:wopt_}
			w_t = \left\{\begin{aligned}& \bar{L}_{q-(t-k)} (Ax_t+Bu_t), \text{ if }k\leq t< k+q\\
				& 0,  \text{ if }t\geq k+q
			\end{aligned}\right.
		\end{equation} 
		with $L_j$, $j \in \{1,\dots,q\}$ given in~\eqref{eq:LK} below.
	\end{algorithmic}
\end{algorithm}
\section{Numerical examples}\label{sec:5}
\par The two numerical examples presented next consider  two scheduler and control pairs. The controller of the first pair is
\begin{equation}\label{eq:uK}
	u_t = K \hat{x}_t, \ \ \hat{x}_t = \left\{\begin{aligned}
		&{x}_t\text{, if }\sigma_t = 1, \\ 
		& (A+BK)	\hat{x}_{t-1}\text{, if }\sigma_t = 0,
	\end{aligned}\right.
\end{equation}
for a given $K$  and the scheduler relies on checking when a weighted norm of $e_t = x_t-\bar{x}_t$ exceeds a threshold, i.e., when
\begin{equation}\label{eq:sched1} \tau_\ell = \min \{k \in \{1,\dots,\bar{h}-1\}| e_{s_\ell+k}^\intercal X e_{s_\ell+k}> \rho \}\}\end{equation}
where $\rho>0$ is a given threshold and $X$ is a given positive semi-definite matrix. Similar schemes appear in, e.g.,~\cite{xu:04},~\cite{molin:10},~\cite{antunes:20}. The control policy of the second pair is
\begin{equation}\label{eq:ut}
		u_t \!=\! K_\gamma \hat{x}_t, \, \hat{x}_t\! =\! \left\{\begin{aligned}
		&\!{x}_t\text{, if }\sigma_t = 1 \\ 
		&\! (A\!+\!BK_\gamma\!+\!L_\gamma(A+BK_\gamma))	\hat{x}_{t-1}\text{, if }\!\sigma_t \!=\! 0.
	\end{aligned}\right.
\end{equation}
so that
$ u_t = u^*_t, \ \ t \in \{s_\ell+1,\dots, s_{\ell+1}-1\} $
with 
$ u^*_t = K(A+BK+L_{\bar{\gamma}}(A+BK))^{t-s_{\bar{\ell}(t)}}x_{s_{\bar{\ell}(t)}}$. To define the scheduler consider a given $\bar{\gamma}$ with $\bar{\gamma}>\gamma_1$. Suppose that the controller is given by~\eqref{eq:uK} with $K=K_{\bar{\gamma}}$ and with $\gamma$ replaced by $\bar{\gamma}$.  The scheduler, adjusted from the one proposed in~\cite{hadi:20}, is defined by
\begin{equation}\label{eq:schhadi}\mu_{\sigma,t}(\mathcal{I}_t) = \left\{\begin{aligned}&1 \text{ if }t = 1, \tau_\ell=\bar{h}  \text{ or }G(\mathcal{I}_t)>0\\
	& 0, \text{ otherwise}.\end{aligned}\right.
\end{equation}
where 

{\footnotesize
$$	\begin{aligned}
		&G(\mathcal{I}_t) = \sum_{k=s_{\bar{\ell}(t)}}^{t-1}(u_{k+1}-{u}^*_{k+1})^\top (R+B^\top F_a(\bar{P}_{\bar{\gamma}})B)(u_{k+1}-{u}^*_{k+1})\\
		& -\!(w_{k}\!-\!\tilde{L}_{\bar{\gamma}}(Ax_{k}\!+\!Bu_{k}))^\top\! ({\bar{\gamma}}^2 I\!-\!\bar{P}_{\bar{\gamma}})^{-1}\!(w_{k}\!-\!\tilde{L}_{\bar{\gamma}}(Ax_{k}\!+\!Bu_{k})).
	\end{aligned}$$
}

\noindent This scheduler-controller pair can be shown to satisfy $\gamma_\pi \leq \bar{\gamma}$ and $r_\pi\leq \frac{1}{h}$~\cite{hadi:20}. 

\subsection{Scalar system}\label{sec:5_1}
\par Suppose that $n=1$, $A=B_2=Q=R=1$. Using~\eqref{eq:iter}, we can compute the  $\gamma_h$ which are given here for $h \in \{1,2,3,4,5\}$:  $\gamma_1=\sqrt{2}$, $\gamma_2=2.0199$, $\gamma_3=2.645$, $\gamma_4=3.276$, $\gamma_5=3.909$. Suppose that $\bar{h}=2$. Then the scheduler must only decide at times $s_\ell+1$, based on $x_{s_\ell}$ and $x_{s_\ell+1}$, if $\sigma_{s_\ell+1}=1$ or $\sigma_{s_\ell+1}=0$. Since $x_{s_\ell+1} = x_{s_\ell}+\mu_{u,s_{\ell}}(x_{s_\ell})+w_k$, we can reparameterize  $\sigma_{s_\ell+1}$ to be a function of $x_{s_\ell}$ and $w_{s_\ell}$ so that we can visualize this decision in $\mathbb{R}^2$.  Consider  $\gamma = (9\gamma_1+\gamma_2)/10=1.4748$. Then ${K}_\gamma=-0.9495$, ${L}_\gamma=8.6463$. Suppose first that the control policy is~\eqref{eq:ut}, and that the scheduler is as in~\eqref{eq:schhadi}, rewritten next based on the numerical values just listed:
$ \sigma_{s_{\ell}+1}= 1\text{ if }\bar{G}(x_{s_{\ell}},w_{s_\ell})>0$, $\sigma_{s_{\ell}}=0$\text{ otherwise}, with
 $ 	 \bar{G}(x,w) = 18.0815(w-0.4366x)^2$
This scheduler only does not trigger transmission in the set of null (Euclidean) measure $(x_t,w_t)\in \{(x,w)|w-0.4366x=0\}$.  With this $\gamma$, Assumption~\ref{as:5} is not met. In fact, if for a given $t=s_\ell$,  $w_t=L_\gamma(Ax_t+Bu_t)=L_\gamma(1+K_\gamma)x_t =0.4366x_t $ then $\bar{G}(x_t,w_t)=0$ and there is no transmission at time $t+1$, but if  $w_t=L_\gamma(Ax_t+Bu_t)+\epsilon$ then $\bar{G}(x_t,w_t) = 18.0815\epsilon^2$ and a transmission will occur irrespective of $\epsilon$.  This can be overcome in two ways. \par
\begin{figure*}
	\includegraphics[width=\textwidth]{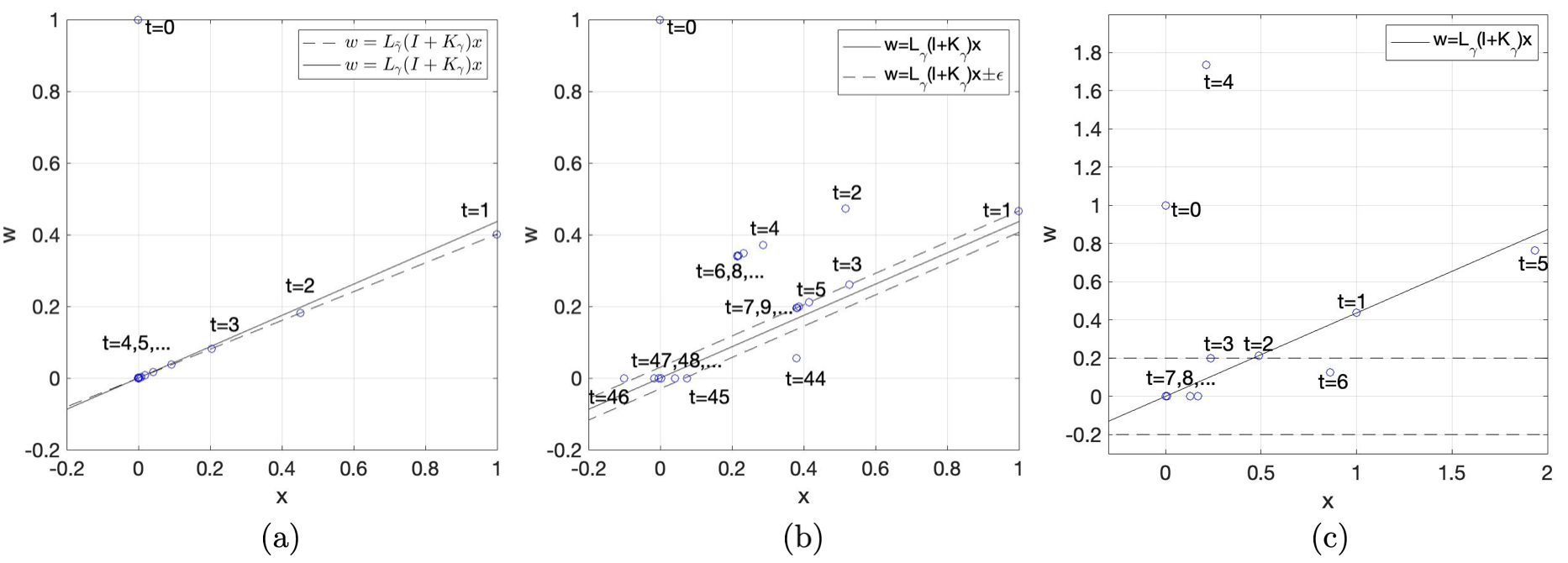}
	\caption{ Trajectories of the scalar system for different schedulers.  At $t=0$, $x_0=0$, $w_0=1$ (arbitrarily picked), leading to $x_1=1$. (a) $\gamma$ in Assumption~\ref{as:5} chosen as $\tilde{\gamma}=1.48$ with  $w_t=L_{\tilde{\gamma}}(A+BK_\gamma)x_t$.  Since transmissions are always triggered $w_t=L_{\gamma_2}(A+BK_\gamma)x_t$ for every $t\geq 2$; (b)  modified scheduler (that meets Assumption~\ref{as:5}) and $\gamma=1.4748$. In this case the disturbance is such that $\|z\|=\gamma \|w\|$.  Between times $t=1$ and $t=43$, Algorithm~\ref{alg:1} probes the system with disturbance~\eqref{eq:w2} with $\epsilon =0$. Since this disturbance is in the $\mathcal{W}$ region (between the dashed lines) it does not lead to transmissions. Thus,  $w_t= L_\gamma(I+K_\gamma)x_t+\epsilon$ is applied, for a sufficiency large $\epsilon$ that still does not lead to transmissions.This happens at time $t=1$ and for odd times until $t=43$ (the pairs $(x_t,w_t)$ are in the border of $\mathcal{W}$ ). At $t=2$ and for even times until $t=43$ there are no transmissions and the disturbance policy is $w_t= L_\gamma(I+K_\gamma)x_t$. At $t=43$, $\eta_t>0$, and Algorithm~\ref{alg:1} applies thereafter  the disturbance policy in Lemma~\ref{lem:key}. Computing $q$ as explained in  Lemma~\ref{lem:key}  leads to $q=1$. After $t=43+q$, the disturbance is always zero and since the control policy is stabilizing the state converges to zero; (c) threshold scheduler. The scheduler triggers when $|w|\leq 0.2$. Thus, for times $t=0$, $t=1$ and $t=2$ the disturbance $w_t = L_\gamma(I+K_\gamma)x_t$ is used. At time $t=3$ since $w_t = L_\gamma(I+K_\gamma)x_t$ would lead to no transmissions the disturbance generator sets the disturbance to $w_3=0.2$ (resulting from $w_t = L_\gamma(I+K_\gamma)x_t+\epsilon$ for some $\epsilon$). This immediately leads to $\eta_t>0$ and this disturbance in Lemma~\ref{lem:key} is applied with $q=3$. After $t=3+q$  the disturbance is zero.}\label{fig:R}
\end{figure*}

 The first is to use the flexibility in picking $\gamma$ allowed by Assumption~\ref{as:5} and test such an assumption with a different $\gamma$, denoted by $\tilde{\gamma}$, such that $\gamma_1<\tilde{\gamma}<\gamma_2$. That is, the scheduler and controller are still the same and pertain to  $\gamma = 1.4748$, but we test Assumption~\ref{as:5} with $\gamma$ replaced by $\tilde{\gamma}=1.48$ and this latter $\tilde{\gamma}$ is the value used for the disturbance policy in Algorithm~\ref{alg:1}. This leads to transmissions being triggered at every time step, and resulting state and disturbances depicted in Figure~\ref{fig:R}(a). At time $t=0$, $x_0=0$, $w_0=1$ (arbitrarily picked), leading to $x_1=1$ at time $1$. This $w_0$ is considered for all the simulations for this scalar example. The scheduler and controller pair ensures $\|z\|^2 = 1.9872  \|w\|^2 < \gamma^2 \|w\|^2$.
\par The second way is to modify the scheduler so that it meets Assumption~\ref{as:5} as already pointed out after Assumption~\ref{as:5}. In this case we do not allow transmissions in the region $\mathcal{W}:=\{(x,w)|w=L_\gamma(I+K_\gamma)x+\epsilon,\epsilon \in  [-\underline{\epsilon},\underline{\epsilon}]\}$, i.e., we change the scheduler to  $\sigma_{s_{\ell}}=0$ if  $\bar{G}(x_{s_{\ell}},w_{s_\ell})\leq 0$ or $w \in \mathcal{W}$, $\sigma_{s_{\ell}}=1$ otherwise. 
Here $\underline{\epsilon}$ can be arbitrarily small; it is set to $\underline{\epsilon}=0.03$. This leads to a disturbance that yields  $\|z\|^2=2.2091 \|w\|^2 > \gamma^2\|w\|^2$. The relevant disturbances and states are plotted in Figure~\ref{fig:R}(b) together with an explanation for such a trajectory.
\par The scheduler and controller pair~\eqref{eq:sched1} is also tested. The same $\gamma=1.4748$ is used, $K$ is set to $K=K_\gamma$,  the threshold to $\rho=0.2$. The policy is simply $\sigma_t=1$ if $|w_{t-1}|>\rho$ when $t=s_{\ell}+1$ for some $\ell$. The scheduler would require a small modification at the intersection point of $ w=L_\gamma(I+K_\gamma)x$ and $w=\pm 0.2$,  but it is  irrelevant when the system trajectory does not pass through these points so this modification is not pursued. In this case the disturbance leads to $\|z\|^2= 2.2726\|w\|^2>  \gamma \|w\|^2$. The state and disturbances are plotted in Figure~\ref{fig:R}(c) together with an explanation for such a trajectory.


\subsection{Third order system}\label{sec:5_2}
\par Suppose that $A      = \begin{bmatrix} 0 & 1 &  0 \\ 0 & 0 & 1 \\ -1  & -2 & -1\end{bmatrix}$, $B = \begin{bmatrix} 0 & 0 &  1\end{bmatrix}^\intercal$, $Q = I_3$, $R = 1$. Then $\gamma_1  = 3.784$, $\gamma_2  = 6.898$, $\gamma_3  = 7.968$, $\gamma_4  = 13.185$,  $\gamma_5  = 15.908$. We set $\gamma=14$ and compare event-triggered control with periodic control with period $4$. As in the previous example using the scheduler-controller pair~\eqref{eq:ut},~\eqref{eq:schhadi}, and picking a slightly different $\gamma$, picked as $\tilde{\gamma}=13.9$ for the disturbance generator leads to all-time transmissions in which the disturbance policy is $w_t=L_{\gamma_2}(A+BK_\gamma)x_t$. We obtain $\|z\|^2=41.49\|w\|^2 \leq \tilde{\gamma}^2\|w\|$.
\par Consider now the controller-scheduler pair~\eqref{eq:sched1}. As in the previous scalar example the scheduler would require a small modification but it is  irrelevant since the system trajectory would not belong to the set requiring modification. The parameter $\bar{\epsilon}$ is set to $\bar{\epsilon}=0.1$. In this case the disturbance leads to $\|z\|^2=197.78\|w\|^2 > \gamma^2\|w\|^2$. 
\section{Concluding remarks}\label{sec:6}
\par We have shown that any event-triggered control strategy, consisting of a scheduler and controller pair,  cannot strictly improve the optimal attenuation bound of periodic control with a smaller or equal average transmission rate. This result was obtained by constructively providing a disturbance policy (Algorithm~\ref{alg:1}) such that for the resulting disturbance either the attenuation bound is larger than or equal to that of the optimal periodic control or the transmission rate is smaller than or equal to that of periodic control. To conclude this, for the proposed disturbance
 policy, we need a technical assumption (Assumption~\ref{as:5}), besides other mild assumptions. Using different proving techniques  it might be possible to obtain the result without  requiring Assumption~\ref{as:5}.

\section{Proofs}\label{sec:7}

\par This appendix provides four auxiliary lemmas. The first  rewrites a special cost using completion of squares and is used by the other three, referred to in the proof of Theorem~\ref{th:1}.
\begin{lemma}\label{lem:completionofsquares}
	Consider~\eqref{eq:sys},~\eqref{eq:output}, $\gamma$ such that $\gamma^2 I-\bar{P}_\gamma$ is invertible, where $\bar{P}_\gamma$ is the unique solution to~\eqref{eq:Pgamma}, and arbitrary  $\tau \in \mathbb{N}$, $\ell \in \mathbb{N}_0$. 
	Then
	{\footnotesize
		\begin{equation}\label{eq:FF3}
			\begin{aligned}
				&\sum_{k=\ell}^{\ell+\tau-1}z_k^\top z_k-\gamma^2 w_k^\top w_k+x_{\ell+\tau}^\top \bar{P}_\gamma x_{\ell+\tau}=x_\ell^\top \bar{P}_\gamma x_\ell+\\
				& \sum_{k=\ell}^{\ell+\tau-1}(u_k-K_\gamma x_k)^\top (R+B^\top F_a(\bar{P}_\gamma)B)(u_k-K_\gamma x_k)-\\
				&\sum_{k=\ell}^{\ell+\tau-1}(w_k-{L}_\gamma(Ax_k+Bu_k))^\top (\gamma^2 I-\bar{P}_\gamma)(w_k-{L}_\gamma(Ax_k+Bu_k))
		\end{aligned}\end{equation}	}
	where $K_\gamma$ and $L_\gamma$ are given by~\eqref{eq:K},~\eqref{eq:Lgamma}. \hfill $\square$
\end{lemma}

\begin{proof2}
	Since~\eqref{eq:sys} is time-invariant it suffices to prove the result for $\ell=0$, which  simplifies the notation. By completion of squares we obtain for every $k \in \{0,\dots,\tau-1\}$
	\begin{equation*}
		\begin{aligned}
			&-\gamma^2 w_{k}^\top w_{k}+x_{k+1}^\top \bar{P}_\gamma x_{k+1}=f_1(x_{k},u_{k},w_{k})+\\
			& (Ax_{k}\!+\!Bu_{k})^\top  \underbrace{(\bar{P}_\gamma + \bar{P}_\gamma(\gamma^2 I-\bar{P}_\gamma)^{-1} \bar{P}_\gamma)}_{=F_a(\bar{P}_\gamma)}(Ax_{k} +Bu_{k})\\
		\end{aligned}
	\end{equation*}
	where
{	\footnotesize
	\begin{equation*}
		\begin{aligned}
			&	f_1(x,u,w) =-(w-L_\gamma(Ax+Bu)  )^\top  (\gamma^2 I-\bar{P}_\gamma) (w-L_\gamma(Ax+Bu)  )
		\end{aligned}
	\end{equation*} }
	Using this equality and again by completion of squares we obtain
	$$\begin{aligned}
		&z_{k}^\top z_{k}  \!-\! \gamma^2w_{k}^\top w_{k}\!+\!x_{k+1}^\top Px_{k+1} \!=\!f_1(x_{k},u_{k},w_{k})\!+\! f_2(x_{k},u_{k})\\
		&+x_{k}^\top \underbrace{(Q+A^\top \tilde{P} A-A^\top  \tilde{P}B(R+B^\top  \tilde{P}B)^{-1}B^\top  \tilde{P}A )}_{\tilde{P}_\gamma}x_{k}.
	\end{aligned}$$
	where $f_2(x,u)=
		(u-K_\gamma x)^\top (R+B^\top F_a(\bar{P}_\gamma)B)(u-K_\gamma x)$
	and $\tilde{P}=F_a(\bar{P}_\gamma)$. Then, using these identities for $k=\tau-1$,
	$$\begin{aligned}
		&\sum_{k=0}^{\tau-1}z_k^\top z_k-\gamma^2 w_k^\top w_k+x_{\tau}^\top \bar{P}_\gamma x_{\tau}=f_2(x_{\tau-1},u_{\tau-1})+\\
		&\sum_{k=0}^{\tau-2}z_k^\top z_k-\gamma^2 w_k^\top w_k\!+\!x_{\tau-1}^\top \bar{P}_\gamma x_{\tau-1}\!+\!f_1(x_{\tau-1},u_{\tau-1},w_{\tau-1})
	\end{aligned}$$
Applying the same procedure for $k=\tau-2$, $k=\tau-3$ until $k=0$ we conclude the desired result. \ \ \ \ \ \ \ \ \ \ \ \ \ \ \ \ \ \ \ \ \ \ \ \ 
\end{proof2}

\begin{lemma}\label{lem:key}
	Consider given $h \in \mathbb{N}$, $\gamma \in \mathbb{R}$, with $\gamma_h< \gamma <\gamma_{h+1}$, so that $\gamma^2 I-M_{h+1}$ has a negative eigenvalue with unitary eigenvector $v_{h,\gamma}$, where $M_h$ is obtained from~\eqref{eq:iter} and $t \in \mathbb{N}$. Suppose that $U_t:=[u_t^\top\, u_{t+1}^\top\, \dots  u^\top_{t+h}]$ are given and consider the following disturbance policy
	\begin{equation}\label{eq:wA}
	w_k=\left\{\begin{aligned}
		&	\bar{L}_\gamma(Ax_k+Bu_k)+\epsilon v_{h,\gamma},\text{ if }k =t, \\
		&	\bar{L}_\gamma(Ax_k+Bu_k),\text{ if }k\in \{t+1,\dots,h\} 
	\end{aligned}\right.
\end{equation}	
for an arbitrary $t \in \mathbb{N}_0$. Then
\begin{equation}\label{eq:FF1}
	\begin{aligned}
	&	\sum_{j=t}^{t+h}z_j^\top z_j-\gamma^2w_j^\top w_j +x_{t+h+1}^\top \bar{P}_\gamma  x_{t+h+1} = x_t^\top \bar{P}_\gamma x_t+\\
	& \ \ \ \ \ \ \ \ \ \ \ \ \ \ \ \ \ \ \ \ \ \ \ \ \ \ a\epsilon^2+b(x_t,U_t)\epsilon+c(x_t,U_t)
	\end{aligned}
\end{equation}
where
	\begin{eqnarray}\nonumber
			 a &\!\!\!=\!\!\!& -v_{h,\gamma}^\top (\gamma^2 I - M_{h+1})v_{h,\gamma} \nonumber \nonumber \\
		    b(x_t,U_t) &\!\!\!=\!\!\!& 2 \sum_{j=t}^{t+h} \phi_{j-t}(x_t,U_t)^\top   Q\bar{A}^{j-t}v_{h,\gamma} \nonumber\\
		    && \ \ \ \ \ \ \ \ \ \ \ +2\phi_{h+1}(x_t,U_t)^\top  \bar{P}_\gamma A^h v_{h,\gamma}\label{eq:FFFF2} \\ 
			c(x_t,U_t) &\!\!\!=\!\!\!& \sum_{j=t}^{t+h} u_j^\intercal Ru_j+\phi_{j-t}(x_t,U_t)^\top Q \phi_{j-t}(x_t,U_t)\nonumber\\
			&& \ \ \ \ \ \ \ \ \ \ \  + \phi_{h+1}(x_t,U_t)^\top \bar{P}_\gamma \phi_{h+1}(x_t,U_t)\nonumber
	\end{eqnarray}
and, letting $\bar{A} = (I+{L}_\gamma)A$, $\bar{B} = (BK+\bar{L}_\gamma B)$,
$$ \phi_j(x_t,U_t) = \bar{A}^jx_t+\sum_{r=0}^{j-1} \bar{A}^{j-1-r}\bar{B}u_{t+r}.$$
Moreover, $a>0$, $c(x_t,U_t)\geq 0$ for every $x_t$, $U_t$. \hfill $\square$
	\end{lemma}
	\begin{proof2}
		We start by noticing that
		$$ x_{t+1}=Ax_t+Bu_t+L_\gamma(Ax_t+Bu_t)+v_{h,\gamma}\epsilon =\bar{A}x_t+\bar{B}u_t+v_{h,\gamma}\epsilon$$
		and for $k>t$, $x_{k+1}=\bar{A}x_k+\bar{B}u_k$. Thus, for $k>t$,
		$$ x_k = \bar{A}^{k-t}v_{h,\gamma}\epsilon+\phi_k(x_t,U_t) .$$
		 The proof then follows by directly replacing this expression on the left hand side of~\eqref{eq:FF1}. The fact that $a>0$ follows from the definition of $v_{h,\gamma}$. The fact that  $c(x_t,U_t)\geq0$ for every $x_t$, $U_t$ follows from~\eqref{eq:FF3} since the left hand side of~\eqref{eq:FF1} when $\epsilon=0$ can be written as a summation of non-negative terms $(u_k-K_\gamma x_k)^\top (R+B^\top F_a(\bar{P}_\gamma)B)(u_k-K_\gamma x_k)$. 
	\end{proof2}
\begin{lemma}\label{lem:2}
	Consider~\eqref{eq:sys} and $\gamma$ such that $\gamma^2 I-\bar{P}_\gamma > 0$ and suppose that Assumption~\ref{as:1} holds. Then
\begin{equation}\label{eq:ine}
\sum_{t=k}^\infty z_t^\top z_t-\gamma^2w_t^\top w_t  \geq x_k^\top G_q x_k,
\end{equation}
	when
	\begin{equation}\label{eq:wopt}
		w_t = \left\{\begin{aligned}& \bar{L}_{q-(t-k)} (Ax_t+Bu_t), \text{ if }k\leq t< k+q\\
			& 0,  \text{ if }t\geq k+q
		\end{aligned}\right.
	\end{equation} 
	where, for $k \in \{1,\dots,q\}$, 
	\begin{equation}\label{eq:LK} \bar{L}_k = (\gamma^2 I-G_{k-1})^{-1}G_{k-1}
	\end{equation} 
	and, for $k \in \{0,\dots,q-1\}$, 
	\begin{equation}\label{eq:Giter}
		G_{k+1}=F_c(F_a(G_k))
	\end{equation}
	with $G_{0} = P_{\text{LQ}}$ where $P_{\text{LQ}}$ is the unique positive definite solution to the algebraic Riccati equation
	$$ P_{\text{LQ}} = A^\top P_{\text{LQ}} A+P_{\text{LQ}}-A^\top P_{\text{LQ}}B(R+B^\top P_{\text{LQ}}B)^{-1}B^\top P_{\text{LQ}}A. $$
	Moreover, for any $x\in \mathbb{R}^n$ and $\beta \in \mathbb{R}_{>0}$, there exists $q \in \mathbb{N}$, denoted by $q=\zeta(x,\beta)$, such that
	\begin{equation}\label{eq:GPgamma}
		\|x^\intercal G_q x-x^\intercal \bar{P}_\gamma x\|<\beta.
	\end{equation} 
	 Such $q$ can be found by running~\eqref{eq:Giter} until~\eqref{eq:GPgamma} is met.   \hfill $\square$
\end{lemma}

\begin{proof2}
	Since~\eqref{eq:sys} is time-invariant it suffices to prove the result for $t=0$, which simplifies the notation. Let 
	$J_{E}(x_r) =\min_{u_t=\mu_{u,t}(\mathcal{J}_t)}\sum_{t=r}^\infty z_t^\top z_t-\gamma^2 w_t^\intercal w_t$
	when $w_t=0$ for every $t\geq r$ in~\eqref{eq:sys}. The standard Linear Quadratic Regulator (LQR) policy is the optimal policy for $u_t$ and leads to the cost $J_{E}(x_r)=x_r^\top P_{LQ} x_r$. 
	 Consider now
	$$\begin{aligned}
		&J_{G,r}(x_0)=\\
		&\min_{u_t=\mu_{u,t}(\mathcal{J}_t)}\max_{w_t=\mu_{w,t}(\mathcal{I}_t)}\sum_{t=0}^{r-1} z_t^\top z_t-\gamma^2w_t^\top w_r+J_{E}(x_r)
	\end{aligned}$$
	From standard arguments for quadratic games~\cite[Ch. 3]{basar:91},   $J_{G,r}(x_0) = x_0^\top G_r x_0$, for optimal disturbance policy~\eqref{eq:wopt} and optimal control policy $u_t=K_{q-(t-k)}x_t$, for $t\in \{k,\dots,k+q\}$, $K_k=-(R+B^\top F_a(G_{k-1})B)^{-1}B^\top F_a(G_{k-1})A.$ This implies~\eqref{eq:ine}.  Moreover, $J_{G,r+1}(x_0) \geq J_{G,r}(x_0)$, and $\lim_{r\rightarrow \infty} G_r = \bar{P}_\gamma$, which implies ~\eqref{eq:GPgamma} is met for some $q$, which can be found with the stated method. 
\end{proof2}

\begin{lemma}\label{lem:final} Consider linear system~\eqref{eq:sys} with performance output~\eqref{eq:output}, an arbitrary $h\in \mathbb{N}$ and suppose that Assumptions~\ref{as:1}-\ref{as:5} hold. Suppose that $w_t$ is generated by Algorithm~\ref{alg:1}.  Then either $\eta_t\leq 0$ for every $s_\ell$, $\ell \in \mathbb{N}$ or $r_\pi(w)\geq \frac{1}{h}$ and $w\in\ell_2$.  
\end{lemma}
\begin{proof2} At times $k=s_\ell$ when~\eqref{eq:A} is  met with $z =x_{s_\ell}$ we have $s_{\ell+1}-s_\ell =j \leq h$ and
$$\sum_{j=s_{\ell}}^{s_{\ell+1}-1}z_j^\top z_j-\gamma^2w_j^\top w_j+ x_{s_{\ell+1}}^\top \bar{P}_\gamma x_{s_{\ell+1}}\geq x_{s_{\ell}}^\top \bar{P}_\gamma x_{s_{\ell}}$$
due to Lemma~\ref{lem:completionofsquares} and the fact that in this case $w_t=L_\gamma(Ax_t+Bu_t)$, $t\in \{s_\ell,\dots,s_{\ell+1}-1\}$. At times $k=s_\ell$ when~\eqref{eq:A} is not met with $z=x_{s_\ell}$ we have $s_{\ell+1}-s_{\ell}>h$ and,
\begin{equation}\label{eq:FF}
	\begin{aligned}
		& \sum_{j=s_{\ell}}^{s_{\ell+1}-1}z_j^\top z_j-\gamma^2w_j^\top w_j+ x_{s_{\ell+1}}^\top \bar{P}_\gamma x_{s_{\ell+1}}\\
		& \ \ \ \ \ \ \ \ \ \ \ = \underbrace{\sum_{j=s_{\ell}}^{s_{\ell}+h-1}z_j^\top z_j-\gamma^2w_j^\top w_j+ x_{s_{\ell}+h}^\top \bar{P}_\gamma x_{s_{\ell}+h}}_{\geq x_{s_{\ell}}^\top \bar{P}_\gamma x_{s_{\ell}} +a\epsilon^2+\epsilon b(x_{s_\ell},U_{s_\ell})+c(x_{s_\ell},U_{s_\ell})  }\\
		& \underbrace{{-x_{s_{\ell}+h}^\top \bar{P}_\gamma x_{s_{\ell}+h}+\sum_{j=s_{\ell}+h}^{s_{\ell+1}-1}z_j^\top z_j-\gamma^2w_j^\top w_j+ x_{s_{\ell+1}}^\top \bar{P}_\gamma x_{s_{\ell+1}}}}_{\geq 0}
	\end{aligned}
\end{equation}
where the first inequality on the right hand side follows from Lemma~\ref{lem:key}, see~\eqref{eq:FF1}, and the second inequality follows from Lemma~\ref{lem:completionofsquares}, see~\eqref{eq:FF3} and the fact that in the interval
 $k\in \{s_{\ell}+h,\dots,s_{\ell+1}-1\}$ the  disturbances in~\eqref{eq:wA2} are given by $w_k=L_\gamma(Ax_k+Bu_k)$. 
Note that due to the choice of $\epsilon$ and the fact that $c(x_{s_\ell},U_{s_\ell})\geq0$, $a>0$, 
\begin{equation}
	a\epsilon^2+\epsilon b(x_{s_\ell},U_{s_\ell})+c(x_{s_\ell},U_{s_\ell})\geq \alpha,
\end{equation} 
with 
$\alpha = v_{h,\gamma}^\top (\gamma^2 I-M)v_{h,\gamma}\bar{\epsilon}^2>0$.
\par Suppose that~\eqref{eq:A} is not met with $\xi=x_{s_\ell}$ $N$ out of $M$ at time steps $s_\ell$, namely $\ell \in \mathcal{N}:=\{i_1,\dots,i_N\}$, $i_p\in  \mathcal{M}:=\{1,\dots,M\}$, and~\eqref{eq:A} is met at times $s_\ell$, $\ell \in \mathcal{M}\setminus \mathcal{N}$. Then, we have
$$ \begin{aligned}
	& \sum_{j=0}^{s_M-1}z_j^\top z_j-\gamma^2w_j^\top w_j+x_{s_M}^\top  \bar{P}_\gamma x_{s_M} 
	\\&= \sum_{j=0}^{s_{1}-1}z_j^\top z_j -\gamma^2 w_j^\top w_j +x_{s_1}^\top  \bar{P}_\gamma x_{s_1} \\ 
	& \ \ \sum_{\ell=1}^{M-1}\!\!(-x_{s_\ell}^\top	\bar{P}_\gamma x_{s_\ell}\!+\!\!\!\!\sum_{j=s_\ell}^{s_{\ell+1}-1}\!\!\!z_j^\top z_j\!-\!\gamma^2w_j^\top w_j\!+\!x_{s_{\ell+1}}^\top  \bar{P}_\gamma x_{s_{\ell+1}})\\
	&=-\gamma^2w_0^\top w_0+\underbrace{\!\sum_{j=1}^{s_{1}-1}z_j^\top z_j -\gamma^2 w_j^\top w_j+ x_{s_1}^\top  \bar{P}_\gamma x_{s_1}}_{\geq \underbrace{x_{1}}_{=w_0}^\top	\bar{P}_\gamma x_{1}\!}+\\ 
	& \sum_{\ell \in \mathcal{N}} \underbrace{(-x_{s_\ell}^\top	\bar{P}_\gamma x_{s_\ell}\!+\!\!\!\!\!\sum_{j=s_\ell}^{s_{\ell+1}-1}\!\!\!z_j^\top z_j-\gamma^2w_j^\top w_j+x_{s_{\ell+1}}^\top  \bar{P}_\gamma x_{s_{\ell+1}})}_{\geq \alpha} +\\
	&\underbrace{\sum_{\ell \in \mathcal{M}\setminus \mathcal{N}}\!\!(-x_{s_\ell}^\top	\bar{P}_\gamma x_{s_\ell}\!+\!\!\!\!\sum_{j=s_\ell}^{s_{\ell+1}-1}\!\!\!z_j^\top z_j\!-\!\gamma^2w_j^\top w_j\!+\!x_{s_{\ell+1}}^\top  \bar{P}_\gamma x_{s_{\ell+1}})}_{\geq 0}\\
	& \geq w_0^\top (-\gamma^2I+ \bar{P}_\gamma )w_0+ \epsilon N
\end{aligned}$$
This implies that~\eqref{eq:A} is not met with $\xi=x_{s_\ell}$ for at most $ N< \delta$ times, where $\delta= \lceil \frac{1}{\epsilon}w_0^\top (\gamma^2I-\bar{P}_\gamma )w_0\rceil$, and where $\lceil a \rceil$ denotes the ceil of a real number $a$. Since $N < \delta$ then $\bar{h}\geq s_{\ell+1}-s_\ell > h$ holds for a finite number of transmission times indexed by $\ell \in \mathbb{N}$. This means that  $s_{\ell+1}-s_\ell \leq h$ holds for an infinite number of transmission times indexed by $\ell \in \mathbb{N}$. This implies that there is $\bar{t}$ such that for $t\geq \bar{t}$  $w_t=L_\gamma(Ax_t+Bu_t)$ and $s_{\ell+1}-s_\ell \leq h$ for $s_\ell\geq \bar{t}$. In turn, this implies that $r_\pi(w)\geq \frac{1}{h}$ as desired. Moreover, it also implies that $w\in \ell_2$. In fact, for  $t\geq \bar{t}$, we have
\begin{equation}\label{eq:FFF1}
	x_{t+1} = \tilde{A}x_t+B\tilde{u}_t, \ \ \tilde{A}=(A+BK_\gamma+L_\gamma(A+BK_\gamma))
\end{equation}
with $\tilde{u}=u_t-K_\gamma(A+BK_\gamma)x_t$ such that $\sum_{j=\bar{t}}^\infty \tilde{u}^\top \tilde{R}\tilde{u}<\infty$ and $\tilde{A}$ is Schur~\cite{basar:91}; this implies that $u\in \ell_2$, that $x_t\in \ell_2$ and thus that $w_t\in \ell_2$. To see  that $\sum_{j=\bar{t}}^\infty \tilde{u}^\top \tilde{R}\tilde{u}<\infty$ holds note that if $\eta_t>0$ is not met for every $t\geq \bar{t}$, from~\eqref{eq:FF3} we conclude that
$$ \eta_t=  \underbrace{\sum_{j=0}^{\bar{t}-1}z_j^\top z_j-\gamma^2w_j^\top w_j+ x_{\bar{t}}^\intercal \bar{P}_\gamma x_{\bar{t}}}_d+\sum_{j=\bar{t}}^{t-1} \tilde{u}^\top \tilde{R}\tilde{u}<0 $$
for every $t\geq \bar{t}$, where $\tilde{R}=R+B^\top F_a(\bar{P}_\gamma)B>0$ and $d$ is a finite constant. Taking the limit as $t\rightarrow\infty$ we conclude $\sum_{j=\bar{t}}^\infty \tilde{u}^\top \tilde{R}\tilde{u}<\infty$.
\end{proof2}

	\bibliography{IEEECL}
	\bibliographystyle{IEEEtran}
\end{document}